\documentclass[a4paper,12pt]{article}

\usepackage{ifthen}
\newboolean{Graphics}

\newcommand{\IfGraphicsOn}[1]{\ifthenelse{\boolean{Graphics}}{#1}{}}

\setboolean{Graphics}{true}
\setboolean{Graphics}{false}

\usepackage{comment}
\usepackage{graphicx}
\usepackage{epstopdf}
\usepackage{color}
\usepackage{enumerate}
\usepackage{url} % to have hypertext-link in the pdf-file
\usepackage{amsmath}
\usepackage{amsfonts}
\usepackage{amssymb,amsthm}
\usepackage{xcolor}
\usepackage{enumerate}
\usepackage{mathrsfs}
\renewcommand\tilde{\widetilde}

\newtheorem{nntheorem}{\bf Theorem}
\newtheorem{nnhypothesis}{\bf Hypothesis}
\newtheorem{nnassumption}{\bf Assumption}
\newtheorem{nndefinition}{\bf Definition}[section]
\newtheorem{nnlemma}[nndefinition]{\bf Lemma}
\newtheorem{nncorollary}[nndefinition]{\bf Corollary}
\newtheorem{nnproposition}[nndefinition]{\bf Proposition}
\newtheorem{nexample}[nndefinition]{\bf Example}
\newtheorem{nnremark}[nndefinition]{\bf Remark}

\newenvironment{thm}
{\begin{nntheorem}\it}
{\end{nntheorem}}

\newenvironment{prop}
{\begin{nnproposition}\it}
{\end{nnproposition}}

\newenvironment{lem}
{\begin{nnlemma}\it}
{\end{nnlemma}}

\newenvironment{dfn}
{\begin{nndefinition}\it}
{\end{nndefinition}}

\newenvironment{assumption}
{\begin{nnassumption}\it}
{\end{nnassumption}}

\newenvironment{note}
{\begin{nnremark}\rm}
{\end{nnremark}}

\newenvironment{nnexample}
{\begin{nexample}\rm}{\end{nexample}}

\DeclareMathOperator{\Lie}{L}

\newcommand{\dx}{\dot{x}}

\newcommand{\mathds}[1]{\mathrm{I\!{#1}}}

\DeclareMathOperator{\Int}{\mathtt{int}}
\DeclareMathOperator{\co}{\mathtt{co}}
\DeclareMathOperator{\dist}{\mathtt{dist}}

\renewcommand{\sup}{\mathtt{sup}\, }
\renewcommand{\max}{\mathtt{max}\, }
\renewcommand{\min}{\mathtt{min}\, }
\renewcommand{\inf}{\mathtt{inf}\, }
\renewcommand{\cos}{\mathtt{cos}\, }
\renewcommand{\sin}{\mathtt{sin}\, }
\DeclareMathOperator{\dom}{\mathtt{dom}}

\newcommand{\Class}{\mathcal{C}}

\newcommand{\Kinf}{\mathcal{K}^\infty}

\newcommand{\R}{\mathbb{R}}

\newcommand{\Ras}{\mathbb{R}_{\geq0}}

\newcommand{\B}{\mathbf{B}}

\newcommand{\U}{\mathbf{U}}

\newcommand{\A}{\mathbf{A}}

\newcommand{\BA}{\mathds{B}}
\newcommand{\HF}{\mathds{K}}
\renewcommand{\H}{\mathds{H}}

\title{\LARGE \bf Combining a backstepping controller with a local stabilizer}
\author{Humberto Stein Shiromoto, Vincent Andrieu, Christophe Prieur
\thanks{Humberto Stein Shiromoto is student from Escola Polit\'{e}nica da Universidade de S\~{a}o Paulo, Avenida Prof. Luciano Gualberto, travessa 3, n 380, CEP 05508-970, S\~{a}o Paulo, SP, Brazil and Politecnico di Torino, Corso Duca degli Abruzzi, 24, 10129 Turin, Italy {\tt humberto.shiromoto@gmail.com}. Vincent Andrieu is with
Universit\'e de Lyon, F-69622, Lyon, France;
Universit\'e Lyon 1, Villeurbanne;
CNRS, UMR 5007, LAGEP.
43 bd du 11 novembre, 69100 Villeurbanne, France.
{\tt https://sites.google.com/site/vincentandrieu/}. Christophe Prieur is with Gipsa-lab,
Department of Automatic Control, 961 rue de la Houille Blanche, BP 46,
38402 Grenoble Cedex, France. {\tt christophe.prieur@gipsa-lab.grenoble-inp.fr}. This work has been initiated during an internship of Humberto Stein Shiromoto at Gipsa-lab, Grenoble.}}
\makeindex

\date{\today}

\begin{document}
 \maketitle

  \begin{abstract}
   We consider nonlinear control systems for which there exist some structural obstacles to the design of classical continuous stabilizing feedback laws. More precisely, it is studied systems for which the backstepping tool for the design of stabilizers can not be applied. On the contrary, it leads to feedback laws such that the origin of the closed-loop system is not globally asymptotically stable, but  a suitable attractor (strictly containing the origin) is practically asymptotically stable. Then, a design method is suggested to build a hybrid feedback law combining a backstepping controller with a locally stabilizing controller. The results are illustrated for a nonlinear system which, due to the structure of the system, does not have {\em a priori} any globally stabilizing backstepping controller.
  \end{abstract}

  \section{Introduction}

  Over the years, research in control of nonlinear dynamical systems has
lead to many different tools to
design (globally) asymptotically stabilizing feedbacks, see
e.g. \cite{FreemanKokotovic:2008,Khalil:1992,Kokotovic:1992}.
Usually these techniques require to impose special structure on the control systems.
Depending on the assumptions made on the model, the designer may use high-gain approaches (as in \cite{GrognardSepulchreBastin:1999}), a backstepping technique (see \cite{FreemanKokotovic:2008,KrsticKokotovicKanellakopoulos_Book_95,PralyAndreaCoron91}) or a forwarding approach (consider e.g., \cite{JankovicSepulchre96,MazencPraly_TAC_96,SepulchreJankovic97}), among others design methods. Unfortunately, in presence of unknown parameters or unstructured dynamics, these classical design methods may fail and some structural obstacles to large regions of attraction may exist. Examples of such systems are the partially linear cascades systems,
considered e.g. in  \cite{BraslavskyMiddleton:ieee96,SaberiKokotovicSussmann:siam90} and \cite{SussmannKokotovic:ieee91}, for which the local stabilization is linear but a perturbation may cause finite escape time, if some parts are not properly controlled. This phenomenon, so-called slow-peaking, has been studied (e.g. in \cite{Sepulchre:ieee:00,SepulchreArcakTeel:Ieee02}) to design
global stabilizers.

For such systems where the classical backstepping techniques can not be applied, the approach presented may solve the problem by combining a backstepping feedback law with a locally stabilizing controller.
More precisely, it is designed a hybrid feedback law to blend both kinds of controllers. The backstepping controller renders a suitable compact set globally attractive, whereas the local one is assumed to have its basin of attraction containing the attractor of the system in closed-loop with the
backstepping controller. The main result can thus be seen as a design techniques of hybrid feedback laws for systems, which {\em a priori} do not have classical nonlinear stabilizing controllers. The use of hybrid stabilizers for systems which do not have continuous stabilizers, is by now classical (see e.g., \cite{HespanhaLiberzonMorse99,MorinSamson00,Prieur:2001}). This approach has been particularly fruitful for control systems that do not satisfy the Brockett's condition \cite{Brockett83} that is a necessary topological condition for the existence of a continuous stabilizing feedback (see in particular
\cite{GoebelHespanhaTeelCaiSanfelice:2004,
GoebelPrieurTeel:2009,HespanhaLiberzonMorse:auto:03,HespanhaMorse99,PrieurAstolfi03}). The considered class of hybrid feedback laws has the advantage to guarantee a robustness property with respect to measurement noise, actuators errors (see \cite{PrieurGoebelTeel:2007} and also \cite{GoebelTeel:2006} for related issues).

Best to our knowledge this is the first work suggesting a design method to adapt the backstepping  technique to a given local controller in the context of hybrid feedback laws. Other works do exist in the context of continuous controllers (e.g., see \cite{PanEzalKrenerKokotovic_TAC_01} where a backstepping controller is blent with an LQ controller, and consider \cite{PrieurAndrieu:2010} where, using control Lyapunov functions, a
globally stabilizing controller is combined with a local optimal controller).
In contrast to these works, for the class of systems considered in this paper, {\em a priori} no continuous stabilizing controller does exist.

This paper is organized as follows. In Section \ref{problem statement}, we introduce precisely the problem under consideration in this paper and the class of controllers that will be used to solve this problem. In Section \ref{main result} the main result is stated, that is the existence of a hybrid feedback law combining a backstepping controller with a local stabilizer. In Section \ref{simu}, the main result is illustrated on an example, and it is designed such a hybrid feedback law for a system for which the classical backstepping approach can not be applied. All technical proofs are collected in Section \ref{sec:proof:proposition}, and Section \ref{conclusion} contains some concluding remarks.

 The proof of some results has been removed due to space limitation.
\begin{comment}
  \emph{Notation}. In this article the Euclidean inner product of two vectors $x$ and $y$ will be denoted by $x\cdot y$, the induced norm will be denoted by $|\cdot|$. The distance from a point $\bar{x}$ to a set $\A$, i.e., $\inf_{y\in Y}|x-y|$ will be denoted by $|\bar{x}|_\A$. The Hausdorff distance of sets $A$ and $B$ will be denoted by $\dist_\text{H}(A,B)$. The closed ball with radius $r$ centered in $x$ will be denoted by $\B_r(x)$. In particular, we denote by $\B$ is the closed unit ball centered in the origin. The interior of a set $X$ will be denoted by $\Int(X)$ and the convex hull $\co(X)$.
  %The interior of a set $S$ will be denoted by $\interior(S)$.
  Regarding the derivative notation, we will write $f'(\bar{x})$ to denote the derivative of a function $f$, at $\bar{x}$. The partial derivative of a function $f$ with respect to $x$ computed at $\bar{x}$ will be denoted by $\partial_x f(\bar{x})$. The Lie derivative of function $V$ with respect to the vector $f$, i.e., $\partial_xV \cdot f$ will be denoted by $\Lie_{f}V$. For a function $f: (x,u)\mapsto f(x,u)$, the notation $f(u)$ stands for $f(x,u)$.
\end{comment}
  \section{Problem statement}\label{problem statement}

  Consider the nonlinear system
  \begin{equation}\label{eq:main system}
   \left\{\begin{array}{rcl}
    \dx_1&=&f_1(x_1,x_2)+h_1(x_1,x_2,u)\\
    \dx_2&=&f_2(x_1,x_2)u+h_2(x_1,x_2,u),\\
   \end{array}\right.
  \end{equation}
  where $(x_1,x_2)\in\R^{n-1}\times\R$, $u\in\U$ is an admissible input. The functions $f_1$, $f_2$, $h_1$ and $h_2$ are locally Lipschitz continuous. Furthermore, the functions satisfy $f_1(0,0)=h_1(0,0,0)=h_2(0,0,0)=0$ and $f_2(x_1,x_2)\neq0$, $\forall (x_1,x_2)\in\R^n$.

  In a more compact notation, we denote system \eqref{eq:main system} by $\dx=f_h(x,u)$. Furthermore, when $h_1\equiv0$ and $h_2\equiv0$ we write $\dx=f(x,u)$.

  \subsection{Assumptions}

  The first assumption concerns the local stabilizability around the origin of system (\ref{eq:main system}). More precisely,
  \begin{assumption}{\em (Local stabilizability)}
  \label{assume:linearization}
  There exist a $\Class^1$ positive definite and proper function $V_\ell:\R^{n}\rightarrow \Ras$, a continuous function $\varphi_\ell:\R^{n}\rightarrow \R$ and a positive constant $v_\ell$ such that,
  $$\partial_{x} V_\ell(x)\cdot  f_h (x, \varphi_\ell (x))< 0\ ,\ \forall x\in \{x:0< V_\ell(x)\leq v_\ell\}.$$
  \end{assumption}

  Note that, when the first order approximation of system (\ref{eq:main system}) is controllable, Assumption \ref{assume:linearization} is trivially satisfied. Indeed, if the couple of matrices $(A,B)$, with $A=\partial_x f_h(0,0)$ and $B=\partial_u f_h(0,0)$ is controllable, then there exist matrices $P>0$ and $K$ such that $V_\ell(x) = x^TPx$ and $\varphi_\ell(x)=Kx$. Thus Assumption \ref{assume:linearization} holds with a sufficiently small positive constant $v_\ell$.

  The second hypothesis provides estimates on terms which prevents using the traditional backstepping method. More precisely, this assumption concerns the global stabilizability of the system
 \begin{equation}\label{eq:x1 subsystem}
   \dot{x}_1=f_1(x_1,x_2)\\
 \end{equation}
 with $x_2$ as an input and bounds of functions $h_1$ and $h_2$. This assumption will be also useful to state a global practical stability property of (\ref{eq:main system}) (see Proposition \ref{th:existence of a controller psi} below).

\begin{assumption}\label{assume:bounded h}
There exist a $\Class^1$ proper and positive definite function $V_1:\R^{n-1}\rightarrow \Ras$, a $\Class^1$ function $\varphi_1:\R^{n-1}\rightarrow \R$ such that $\varphi_1(0)=0$, a locally Lipschitz $\Kinf$ function $\alpha:\Ras\rightarrow\Ras$,
a continuous function $\Psi:\R^n\rightarrow \R$ and two positive constants $\varepsilon<1$ and $M$ such that the following properties hold.
\begin{enumerate}
    \item {\em (Stabilizing controller $\varphi_1$ for (\ref{eq:x1 subsystem}))}
  $\forall x_1\in\R^{n-1}$,
$$
\partial_{x_1} V_1 (x_1)\cdot f_1 (x_1, \varphi_1 (x_1))\leq -\alpha(V_1(x_1)).
  $$
% is a proper positive definite function;
 \item {\em (Estimation on $h_1$)} $\forall(x_1,x_2,u)\in\R^{n-1}\times \R\times \R$,
  \\[0.5em]$L_{h_1} V_1(x_1,\varphi_1(x_1),u)\leq(1-\varepsilon )\alpha(V_1(x_1))$\\[0.5em]
  \null\hfill$+\varepsilon \alpha(M),$\hfill\refstepcounter{equation}$(\theequation)$\label{estimation1onh1}\\[0.5em]
  \null\hfill$| h_1(x_1,x_2,u)|\leq \Psi(x_1,x_2)$\hfill\refstepcounter{equation}$(\theequation)$\label{estimation2onh1}

 \item {\em (Estimation on $\partial_{x_2}h_1$)}
 $\forall(x_1,x_2,u)\in\R^{n-1}\times \R\times \R$,
 \begin{equation}
 \label{estimation1onpartialh1}
 |\partial_{x_2} h_1(x_1,x_2,u)| \leq \Psi(x_1,x_2).
 \end{equation}

\item  {\em (Estimation on $h_2$)}
$\forall(x_1,x_2,u)\in\R^{n-1}\times \R\times \R$,
 \begin{equation}
 \label{estimation1onh2}
| h_2(x_1,x_2,u)|\leq \Psi(x_1,x_2).
 \end{equation}

\end{enumerate}
\end{assumption}

As we will see in this work, it is not necessary that $\varphi_1$ be $\Class^1$ in a neighborhood of the origin because, in such a region, we use the  local controller $\varphi_\ell$.

Before introducing the third assumption, let us denote $\A$ the subset of $\R^n$ defined by
\begin{equation}
\label{def:attractor}
\mathbf{A}=\{(x_1, x_2) \in\R^{n}:V_1(x_1) \leq
M,\ x_2=\varphi_1(x_1)\}.
\end{equation}
Note that since, by Assumption \ref{assume:bounded h}, the function $V_1$ is proper, this set is compact.
Moreover, it will be proven below (see Proposition \ref{th:existence of a controller psi}) that with the other items of Assumption \ref{assume:bounded h} a controller to (\ref{eq:main system}) can be designed such that $\A$ is globally practically stable to the system in closed-loop with this controller.

The last assumption describes that $\A$ is included in the basin of attraction of the controller $\varphi_\ell$.
%More precisely, it is assumed the following assumption.
\begin{assumption}{\em (Inclusion  assumption)}
\label{assume:Covering}
\begin{equation}
\label{last:condition:assume:2}
\max_{x\in \mathbf{A}} V_\ell (x)< v_\ell\ .
\end{equation}
\end{assumption}

The problem under consideration in this paper is the design of a controller such that the origin is globally asymptotically stable for \eqref{eq:main system}. Due to the presence of the functions $h_1$ and $h_2$ and their dependence with respect to $u$, a classical backstepping can not be achieved to compute a global stabilizer.%
\footnote{\label{note:1} More precisely, following the classical basckstepping approach, let us assume that item 1 of Assumption \ref{assume:bounded h} holds and let us consider the Lyapunov function candidate $V(x_1,x_2)=V_1(x_1) + \frac{1}{2}(x_2-\varphi_1(x_1))^2$. We compute along the solutions of \eqref{eq:main system}, for all $(x_1,x_2,u)$ in $\R^{n-1}\times \R\times \R$,
$$
\begin{array}{rcl}
\dot V &\leq &-\alpha(V_1(x_1))+ \left[x_2 -\varphi_1(x_1)\right]\left[f_2(x_1,x_2)u +h_2(x_1,x_2,u)\right.\\
&& \left.-\frac{\partial \varphi_1}{\partial x_1}(x_1)\cdot (f_1 (x_1,x_2) + h_1 (x_1,x_2 ,u))
\right.\\
&&\left.
+\frac{\partial V_1}{\partial x_1}(x_1)\cdot \int_0 ^1 f _1(x_1,s x_2 - (1-s) \varphi_1 (x_1) ) ds\right]
\\
&& 
+ \frac{\partial V_1}{\partial x_1}(x_1)\cdot h_1 (x_1,x_2,u) \ .
\end{array}
$$
And thus to get an term $(x_2- \varphi_1(x_1))^2$ in the right-hand side of this inequality, it is natural to look for a control
$u=u(x_1,x_2)$ satisfying the following identity, for all $(x_1,x_2)$ in $\R^{n-1}\times \R$,
$$\begin{array}{c}
f_2(x_1,x_2)u +h_2(x_1,x_2,u) -\frac{\partial \varphi_1 }{\partial x_1}(x_1)\cdot (f_1 (x_1,x_2) + h_1 (x_1,x_2 ,u))
\\
+\frac{\partial V_1}{\partial x_1}(x_1)\cdot \int_0 ^1 f _1(x_1,s x_2 - (1-s) \varphi_1 (x_1) ) ds = - k\left(x_2 -\varphi_1(x_1)\right)
\end{array}
$$
for some positive value $k$. However this equation is implicit in the variable $u$ due to dependance of $h_1$ and of $h_2$
with respect to $u$. Therefore it seems to us that the classical backstepping cannot be achieved to compute a stabilizer for \eqref{eq:main system}.
}

However we succeed to design a controller rendering a compact set globally asymptotically stable to (\ref{eq:main system}) in closed-loop. Then a natural approach is to combine this controller with a local feedback law given by Assumption \ref{assume:linearization}. Global asymptotical stabilization of the origin of $\R^n$ can be achieved by considering a hybrid controller which blends the different controllers according to each basin of attraction. The strategy is similar to that one developed in \cite{Prieur:2001}, namely, we divide the continuous state space in two open sets introducing a region with hysteresis. This asks to make precise the class of controllers under consideration in this paper.

\subsection{Class of controllers}

\begin{dfn}
     A \emph{hybrid feedback law} to \eqref{eq:main system}, denoted by $\HF$, consists of
     \begin{itemize}
      \item a totally ordered countable set $Q$;
      \item for each $q\in Q$,
      \begin{itemize}
       \item closed sets $C_q\subset\R^n$ and $D_q\subset\R^n$ such that $C_q\cup D_q=\R^n$;
       \item a continuous function $\varphi_q:C_q\to\R$;
       \item an outer semi-continuous%
       \footnote{a set-valued mapping $F: \R^m\rightrightarrows \R^n$ is said to be outer 
       semicontinuous if each sequence $(x_i,f_i)$ in $\R^m\times \R^n$ that satisfies $f_i\in F(x_i)$ for each $i$, and 
       converges to a point $(x,f)$ in $\R^m\times \R^n$ has the property that $f\in F(x)$.}, %
       %end footnote
and  locally bounded%
\footnote{a set-valued mapping $F: \R^m\rightrightarrows \R^n$ is said to be locally bounded if, for each compact set $K_1\subset \R^n$, there exists a compact set $K_2\subset \R^n$ such that $F(K_1):=\bigcup _{x\in K_1} \subset K_2$. The boundedness is said to be uniform with respect to a parameter if the set $K_2$ can be selected uniformly with respect to this parameter.},
% end footnote
 uniformly in $q$, set-valued mapping 
       $G_q:D_q\rightrightarrows Q$ with non-empty images,
      \end{itemize}
     \end{itemize}
     such that the family $\{C_q\}_{q\in Q}$ is locally finite and covers $\R^n$. 
    \end{dfn}

% We assume that, for each $q\in Q$, $\varphi_q$ is a measurable, essentially bounded function satisfying $x\mapsto%\varphi_q(x)\in\R$, a.e. $\ x\in C_q$. We also consider the Assumptions A$_q$1)-A$_q$3), A4)-A6) from \cite{PrieurGoebelTeel:2007}, that is we assume that the sets $C_Q$ and $D_q$ are closed for each $q$ in $Q$

  System \eqref{eq:main system} in closed loop with $\HF$ lies in the class of hybrid systems as considered in e.g., \cite{Aubinetal:2002}. It is defined as the hybrid system
  \begin{equation}\label{eq:hybrid system}
    \H:\left\{\begin{array}{rcll}
      \dot{x}&=&f_h(x,\varphi_q(x)),&x\in C_q\\
      q^+&\in&G_q(x),&x\in D_q.
    \end{array}\right.
  \end{equation}
  Note that the state space of $\H$ is $\R^n\times Q$.

  \begin{dfn}
   A \emph{hybrid time domain} $S\subset\Ras\times\mathbb{N}$, is the union of finitely of infinitely many time intervals $[t_j,t_{j+1}]\times\{j\}$, where the sequence $\{t_j\}_{j\geq0}$ in nondecreasing, with the last interval, if it exists, possibly of the form $[t,T)$ with $T$ finite or $T=\infty$.
  \end{dfn}

  \begin{dfn}
   A \emph{solution} to $\H$ with initial condition $(x(0,0),q(0,0))=(x_0,q_0)$ consists of
   \begin{itemize}
    \item A hybrid time domain $S\neq\emptyset$;

    \item A function $x:S\to\R^n$, where $t\mapsto x(t,j)$ is absolutely continuous, for a fixed $j$, and constant in $j$ for a fixed $t$ over $(t,j)\in S$;

    \item A function $q:S\to Q$ such that $q(t,j)$ is constant in $t$, for a fixed $j$ over $(t,j)\in S$;

   \end{itemize}
   meeting the conditions
    \begin{enumerate}[{S}$_1$)]
     \item $x(0,0)\in C_{q(0,0)}\cup D_{q(0,0)}$;

     \item $\forall j\in\mathbb{N}$ and $\ae\ t$ such that $(t,j)\in S$,
	   $$\dot q(t,j)=0, \ \dx(t,j)\in F_{q(t,j)}(x(t,j)),\ x(t,j)\in C_{q(t,j)};$$

      \item $\forall (t,j)\in S$ such that $(t,j+1)\in S$
	   $$\begin{array}{c}
	   x(t,j+1)=x(t,j), \ q(t,j+1)\in G_{q(t,j)}(x(t,j)),\\ x(t,j)\in D_{q(t,j)}.\end{array}$$
    \end{enumerate}
  \end{dfn}

  From now on, we will refer to the domain of a solution $(x,q)$ to $\H$ as $\dom(x,q)$. A solution $(x,q)$ to $\H$ is called {\em maximal} if it cannot be extended, i.e., does not exists any solution defined on a larger domain of definition and equal to $(x,q)$ on $\dom(x,q)$. A solution is {\em complete} if its domain is unbounded.

  During flows, $x$ evolves according to the differential equation $\dx=f_h(x,\varphi_q(x))$, $x\in C_q$ while $q$ remains constant. During jumps, $q$ evolves according to the difference inclusion $q^+\in G_q(x)$, $x\in D_q$ while $x$ remains constant.

  \begin{note}
   Note that a sufficient condition for the existence of a hybrid stabilizer for \eqref{eq:main system} is the global asymptotic controllability (see \cite{PrieurGoebelTeel:2007}, Theorem 3.7 for more details).

   Together with locally Lipschitz continuity assumption, we consider the Filippov regularization of \eqref{eq:main system} which assures existence, uniqueness and bounded dependence on the initial condition for solutions of $\H$. Moreover, $\H$ is robust and its solution behaves as follows: it is either  complete or blows in a finite hybrid domain time or eventually jumps out of $C_q\cup D_q$, $q\in Q$. For further information, see \cite{AubinCellina:1984}, \cite{BacciottiRosier:2010}, \cite{Filippov:1988}, \cite{GoebelSanfeliceTeel:2009} and \cite{GoebelTeel:2006}.
  \end{note}

  We can now define the notion of stability needed to design the controller for the hybrid closed loop system.
\begin{dfn}\hfill
  \begin{itemize}
    \item A set $\A\subset \R^n$ is {\em stable} for $\H$ if $\forall\varepsilon>0$, $\exists\delta>0$ such that any solution $(x,q)$ to \eqref{eq:hybrid system} with $|x_0|_\A\leq\delta$ satisfies $|x(t,j)|_\A\leq\varepsilon$, for all $(t,j)\in\dom(x,q)$;

   \item A set $\A\subset \R^n$ is {\em attractive} for $\H$ if there exists $\delta>0$ such that
   \begin{itemize}
    \item for all $(\bar{x},\bar{q})\in\R^n\times Q$ with $|\bar{x}|_\A\leq\delta$ there exists a solution to $\H$ with $(x,q)(0,0)=(\bar{x},\bar{q})$;
    \item for any maximal solution $(x,q)$ to $\H$ with $|x(0,0)|_\A\leq\delta$ we have $|x(t,j)|_\A\to0$ as $t\to\sup_t(\dom(x,q))$.
   \end{itemize}

   \item The set $\A\subset \R^n$ is {\em asymptotically stable} if it is stable and attractive;

   \item The {\em basin of attraction}, denoted by $\BA_\H(\A)$, is the set of all $\bar{x}\in \R^n$ such that for all $\bar{q}\in Q$, there exists a solution to $\H$ with $x(0,0)=\bar{x}$, $q(0,0)=\bar{q}$ and any such solution that is also maximal satisfies $|x(t,j)|_\A\to0$ as $t\to\sup_t\dom(x,q)$;

  \item The set $\A\subset \R^n$ is {\em globally asymptotically stable} if  $\BA_\H(\A)=\R^n$.
  \end{itemize}
\end{dfn}

  \section{Main result}\label{main result}

  Let us denote the unit closed ball in $\R^n$ by $\B$. Before stating our main result, let us first solve a preliminary design problem by adapting the backstepping technique:

\begin{prop}\label{th:existence of a controller psi}
Under Assumption \ref{assume:bounded h}, the set $\A$ defined by (\ref{def:attractor}) is globally practically stabilizable, i.e. for each $a>0$ there exists a continuous controller $\varphi_g$
such that the set $\A+a\B$ contains a set that is globally asymptotically stable for system (\ref{eq:main system}) in closed-loop with $\varphi_g$.
\end{prop}
\begin{comment}
The proof of Proposition \ref{th:existence of a controller psi} has been removed due to space limitation.
However, this proof is constructive and given a positive real number $a>0$ the controller $u=\varphi_g(x)$ 
guaranteeing global asymptotic stability of the set  $\mathcal{A}+a\mathcal{B}$ is given as
 \\[0.5em]$\varphi_g(x_1,x_2,a)=\frac{1}{f_2(x_1,x_2)}\left[\frac{\tilde{u}}{k}+\partial_{x_1}\varphi_1(x_1)\cdot f_1(x_1,x_2)\right.$\\[0.5em]
 \null\hfill$\left.+\frac{1}{k}\partial_{x_1}V_1(x_1)\cdot\int_0^1\partial_{x_2}f_1(x_1,\eta_{x_1,x_2}(s))\,ds\right],$\hfill\refstepcounter{equation}$(\theequation)$\label{eq:global controller} \\[0.5em]
  where
\\[0.5em]\null\hfill$\tilde{u}=-(x_2-\varphi_1(x_1))c[1+\frac{1}{4}\Delta(x_1,x_2)^2],$\hfill\refstepcounter{equation}$(\theequation)$\label{eq:tildeu}
\\[0.5em]\null\hfill$
k=2\frac{M+a}{a^2}
$\hfill\refstepcounter{equation}$(\theequation)$\label{def:k}\\[0.5em] and
\\[0.5em]$\Delta(x_1,x_2)=|\partial_{x_1}V_1(x_1)|\int_0^1\Psi(x_1,\eta_{x_1,x_2}(s))\,ds$\\[0.5em]
\null\hfill$+\Psi(x_1,x_2)k(1+|\partial_{x_1}\varphi_1(x_1)|)$\hfill\refstepcounter{equation}$(\theequation)$\label{eq:Delta}\\[0.5em]
\end{comment}

We are now in position to state our main result.
\begin{thm}\label{main:result}
 Let $v_\ell$ and $\tilde{v}_\ell$ be two positive constants satisfying $0<\tilde{v}_\ell<v_\ell$. 
 Under Assumptions  \ref{assume:linearization}, \ref{assume:bounded h} and \ref{assume:Covering} there exists $a>0$ such that the hybrid controller $\HF$ defined by
 $Q=\{1,2\}$, subsets
 $$
  \begin{array}{rcl}
      C_1&=&\{(x_1,x_2)\in \R^{n-1}\times \R:V_\ell(x_1,x_2)\leq v_\ell\},\\
      C_2&=&\{(x_1,x_2)\in \R^{n-1}\times \R:V_\ell(x_1,x_2)\geq \tilde{v}_\ell\},\\
      D_q&=&\overline{(\R^{n-1}\times \R)\setminus C_q}, \; \forall q=1,2,
     \end{array}
 $$
 controllers $C_1\ni(x_1,x_2)\mapsto\varphi_1(x_1,x_2)=\varphi_\ell(x_1,x_2)\in\R$ and $C_2\ni(x_1,x_2)\mapsto\varphi_2(x_1,x_2)=\varphi_g(x_1,x_2,a)\in\R$ and set-valued mapping $D_q\ni(x_1,x_2)\mapsto G_q(x_1,x_2)=\{3-q\}$, $q\in Q$, renders the origin globally asymptotically stable for \eqref{eq:main system} in closed-loop with $\HF$.

\end{thm}

Let us emphasize that this result is more than an existence result since its proof allows to design a suitable hybrid feedback law. Let us sketch the proof of Theorem \ref{main:result}.
First, we use Assumption \ref{assume:bounded h}, and Proposition \ref{th:existence of a controller psi} is applied to design a controller, denoted $\varphi_g$, such that the set $\A$ is globally practically stable for the system \eqref{eq:main system} in closed-loop with $\varphi_g$. Using Assumptions \ref{assume:linearization} and  \ref{assume:Covering}, this set is shown to be included in the basin of attraction of the system \eqref{eq:main system} in closed-loop with $\varphi_\ell$. Then we design a hybrid feedback law based on an hysteresis of both controllers $\varphi_\ell$ and $\varphi_g$ on appropriate sets. This latter construction is adapted from other works like \cite{GoebelSanfeliceTeel:2009} or \cite{Prieur:2001}. The complete proof of Theorem \ref{main:result} is in Section \ref{sec:proof:proposition} below.

\section{Illustration}
\label{simu}

  Before applying the main result of this paper, let us first consider the following example in $\R^2$
 \begin{equation}\label{eq:example:preliminar}
  \left\{\begin{array}{rcl}
   \dot{x}_1&=&x_1+\theta x_1^2+x_2\\
   \dot{x}_2&=&u
  \end{array}\right.\ ,
 \end{equation}
 where $\theta$ is a positive constant.

This system is in backstepping form and many references on how to design a global stabilizer are presented in the literature, for instance, the reader may see \cite{FreemanKokotovic:2008,Khalil:1992}, and \cite{Kokotovic:1992}.
Following this approach, in a first step, we consider the two smooth functions $\varphi_1(x_1)=-(1+c_1)x_1-\theta x_1^2$ and $V_1(x_1) = \frac{1}{2} x_1^2$ where $c_1$ is a positive constant. It can be checked that this function is such that, for all $x_1$ in $\R$,
\begin{equation}\label{ex:varphi1}
\partial_{x_1}V_1(x_1)\left[x_1+\theta x_1^2+\varphi_1(x_1)\right]=-2c_1 V_1(x_1)\ .
\end{equation}

This gives the control law, for all $(x_1,x_2)$ in $\R^2$,
$$\begin{array}{rcl}
   \varphi_b(x_1, x_2)&=&-(1+c_1 + 2 \theta x_1)(x_1 + \theta x_1^2 + x_2)\\
   &&-x_1-c_2\left(x_2+(1+c_1)x_1+ \theta  x_1^2\right)
  \end{array}
  $$
which is such that along the solutions of (\ref{eq:example:preliminar}),
\\[0.5em]$
\dot V_b(x_1, x_2) = -c_1 x_1^2 - c_2\left(x_2 + (1+c_1)x_1+ \theta  x_1^2\right)^2
$\\[0.5em]
where $
V_b(x_1, x_2) = V_1(x_1) +\frac{1}{2} (x_2 +(1+c_1)x_1+ \theta x_1^2)^2$.

However the backstepping technique cannot be applied to the following system:
\begin{equation}\label{eq:example:1}
  \left\{\begin{array}{rcl}
   \dot{x}_1&=&x_1+x_2+\theta [x_1^2+(1+x_1)\sin(u)]\\
   \dot{x}_2&=&u
  \end{array}\right.
 \end{equation}
due to the presence of the term $(1+x_1)\sin(u)$ in the time-derivative of $x_1$ (recall the discussion in Footnote \ref{note:1}). Therefore, it is necessary to revise the controller design for (\ref{eq:main system}) and to apply Theorem \ref{main:result}. With obvious definitions of the functions $f_1$, $f_2$, $h_1$ and $h_2$, system (\ref{eq:example:1}) may be rewritten as system (\ref{eq:main system}) and system (\ref{eq:example:preliminar}) may be rewritten as $\dot x= f(x,u)$. There exists $\theta>0$ sufficiently small such that we may apply Theorem \ref{main:result}. Indeed we have the following result.
\begin{lem}\label{LEM:TECH}
Let $\theta$ be a positive constant. If $\theta$ is sufficiently small, then Assumptions \ref{assume:linearization}, \ref{assume:bounded h}, and \ref{assume:Covering} hold for system \eqref{eq:example:1}.
\end{lem}

The proof has been removed due to space limitation.

Combining this result with Theorem \ref{main:result}, we may design a hybrid feedback law $\HF$ such that
the origin is globally asymptotically stable to (\ref{eq:example:1}) in closed-loop with $\HF$.

Let us consider the following parameters $\theta=10^{-3}$, $\rho=2$, $c_1=\tfrac{(2+\rho)\theta}{2}+1=1.0020$, $a=10$ and $c=10$. Item 1 of Assumption \ref{assume:bounded h} is satisfied with $\alpha(s)=2c_1s$, $\forall s\geq0$. Item 2 is satisfied with positive constants $\varepsilon = 1-\theta\frac{2+\rho}{2c_1}=0.998$ and $M = \frac{\theta}{2\rho(2c_1-\theta(2+\rho))}=1.25\times10^{-4}$ . Items 3 and 4 are satisfied with $\Psi(x_1,x_2)=\theta(1+|x_1|)$.

Since the pair of matrices $(A,B)=(\partial_xf_h(0,0),\partial_uf_h(0,0))$ is controllable, Assumption \ref{assume:linearization} holds with $\varphi_\ell(x)=k_1x_1+k_2x_2$, where $k_1=-5-\theta$ and $k_2=-3+3\theta+\theta^2$, $V_\ell(x)=\frac{1}{2}(x_1-\theta x_2)^2 + \frac{1}{2}(2x_1+(1-2\theta) x_2)^2$ and $v_\ell = \left(\frac{2}{\theta p(\theta)}\right)^2=0.1042$. Moreover, in the set defined by
$$\A=\{(x_1,x_2)\in\R^2:|x_1|\leq\sqrt{2.5\times10^{-4}},\ x_2=\varphi_1(x_1)\}.$$
we may check that
$$\max_{x\in\A}V_\ell(x_1,x_2)=0.0001<v_\ell,$$
and thus Assumption \ref{assume:Covering} holds. Following Proposition \ref{th:existence of a controller psi} and Theorem \ref{main:result}, we may define a hybrid controller. More precisely, computing $k=2\frac{M+a}{a^2}=0.2$, we define the global controller
\\[0.5em]
$\varphi_g(x_1, x_2)=\frac{\tilde u}{k}-(1+c_1+2\theta x_1)(x_1+\theta x_1^2+x_2)+\frac{x_1}{2k},$\\[0.5em]
where $\tilde{u}=(x_1-\varphi_1(x_1))\left[-c-\frac{c}{4}\Delta(x_1,x_2)^2\right]$ and
\\[0.5em]$\Delta(x_1,x_2)=|x_1|\theta(1+|x_1|)+\theta(1+|x_1|)$\\[0.5em]
\null\hfill$\cdot k(1+|(1+c_1)x_1+\theta x_1^2|)$

Then, letting $\tilde{v}_\ell=0.05$, the origin is globally asymptotically stable for \eqref{eq:example:1} in closed-loop with the hybrid controller $\HF$ defined as in  Theorem \ref{main:result}.

\begin{figure}[htb!]
 \centering
\includegraphics[width=\textwidth]{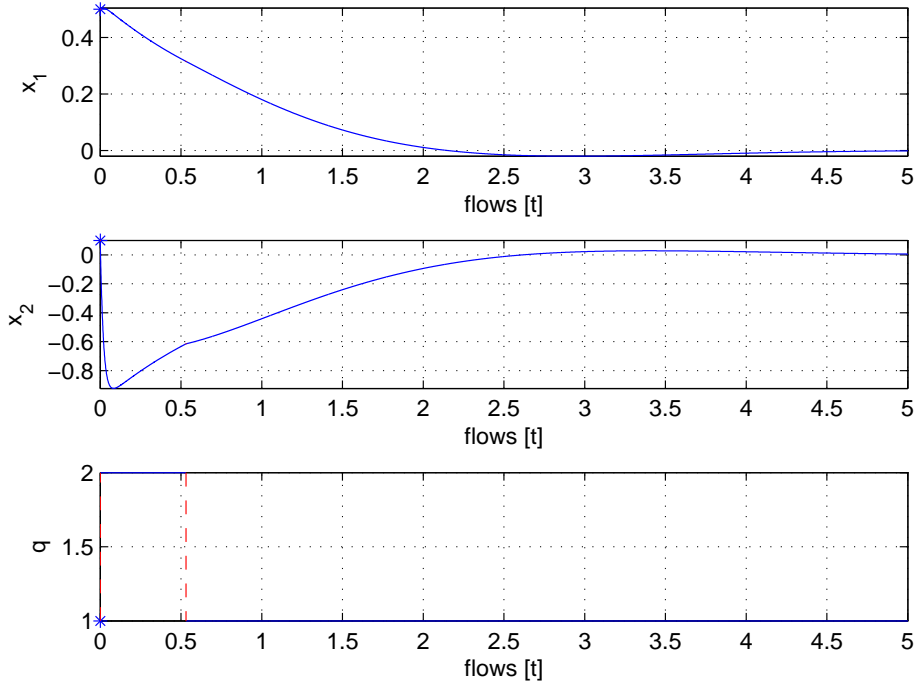}
 % x_q.pdf: 448x337 pixel, 72dpi, 15.80x11.89 cm, bb=0 0 448 337
 \caption{At top, time evolution of $x_1$, at middle, time evolution of $x_2$ and, at bottom, time evolution of $q$.}
 \label{fig:xq}
\end{figure}

Let us check this property on numerical simulations. To do that, we consider the initial condition $x_1(0,0)=0.5$, $x_2(0,0)=0.1$ and $q(0,0)=1$. See Fig. \ref{fig:xq} for the time evolution of the $x_1$, $x_2$ and $q$ components of the solution of  \eqref{eq:example:1} in closed-loop with $\HF$.
First the system (\ref{eq:example:1}) is in closed-loop with the controller $\varphi_g$ (for continuous time between $0$ and  $0.5314$). Then the system (\ref{eq:example:1}) is in closed-loop with the controller $\varphi_\ell$, and the solution converges to the origin. 
\begin{comment}
Fig. \ref{fig:Vlq} gives the value of the function $V_\ell$ along the continuous time. This value is used in the proof of Theorem \ref{main:result} to define the hybrid controller $\HF$.
\end{comment}

\begin{comment}
\begin{figure}[htb!]
 \centering
 \includegraphics[width=\textwidth]{simus(Vlq)_09-26}
 % V.pdf: 448x337 pixel, 72dpi, 15.80x11.89 cm, bb=0 0 448 337
 \caption{Time evolution of $V_\ell$.}
 \label{fig:Vlq}
\end{figure}
\end{comment}
\section{Proof of Theorem \ref{main:result}}
\label{sec:proof:proposition}
%\begin{comment}
\subsection{Proof of Proposition \ref{th:existence of a controller psi}}
\begin{proof}
Let $a$ be a positive value. We wish to show that there exists a continuous controller $\varphi_g$ such that $\A+a\B$ contains a set that is globally and asymptotically stable.

First of all, note that if we introduce the function $r_1(x_1,x_2, u)=f_1(x_1,x_2)+h_1(x_1,x_2,u)$, we get with Item 1 and  Item 2 of Assumption \ref{assume:bounded h} that along the solutions of (\ref{eq:main system}), we have for all $(x_1, x_2)$ in $\R^n$ and $u$ in $\R$,
\begin{equation}\label{eq:dotV1}
  \begin{array}{rcl}
   \dot V_1(x_1)&\leq&\varepsilon[ \alpha(M)-\alpha(V_1(x_1))]\\
   \multicolumn{3}{r}{+\partial_{x_1} V_1(x_1)\cdot [r_1(x_1,x_2, u)- r_1(x_1,\varphi_1(x_1), u)]}
  \end{array}
\end{equation}

Moreover, with the $\Class^1$ function $\eta_{x_1,x_2}(s)=sx_2+(1-s)\varphi_1(x_1)$, it yields
$$\partial_s r_1(x_1,\eta_{x_1,x_2}(s),u)=\partial_{x_2} r_1 (x_1,\eta_{x_1,x_2}(s),u) (x_2-\varphi_1(x_1))\ ,$$
which implies
$$\begin{array}{rcl}
   r_1(x_1,x_2,u)-r_1(x_1,\varphi_1(x_1),u)&=&
   (x_2-\varphi_1(x_1))\\
   \multicolumn{3}{r}{\int_0^1\partial_{x_2} r_1(x_1,\eta_{x_1,x_2}(s),u)\,ds.}
  \end{array}$$
Hence, Equation \eqref{eq:dotV1} becomes,
$$\begin{array}{rcl}
   \dot V_1(x_1)&\leq&\varepsilon[ \alpha(M)-\alpha(V_1(x_1))]\\
   \multicolumn{3}{r}{+ (x_2-\varphi_1(x_1)) \partial_{x_1} V_1(x_1)\cdot \int_0^1\partial_{x_2} r_1(x_1,\eta_{x_1,x_2}(s),u)\,ds.}
  \end{array}
$$

Let $V(x)=V_1(x_1)+\frac{k}{2}(x_2 -\varphi_1(x_1))^2$ for all $(x_1,x_2)$ in $\R^{n}$ with $k=2\frac{M+a}{a^2}$.
%Note that this function is such that if $V(x)\leq M+a$ then $x$ is in $\mathcal A+a\mathcal B$.
Let $a'$ be a positive value such that $V_1(x_1)\leq a'$ implies $x_1\in \{x_1':V_1(x_1')\leq a' \}+a \B$, in other words, $a'$ is such that
$$
V_1(x_1)\leq a'
\Rightarrow \exists x_1' \mbox{ s.t. } V_1 (x_1 ') \leq a' \mbox{ and } |x_1 -x_1' |\leq a
\ .
$$
Such positive value $a'$ exists since $V_1$ is assumed to be a proper function. Let $\tilde a =\min\{ a,a'\}$.  With these definitions of $k$ and $a'$, we get
\begin{equation}
\label{inclus}
\{ x:V(x) \leq M+\tilde a\}\subset \A+ a \B
\end{equation}
%$V_1(x_1)\leq M+a$ and $\frac{M+a}{2a^2}(x_2 -\varphi_1(x_1))^2\leq M+a$
Consider now the control $\varphi_g$ defined for all $\tilde u$ in $\R$ as in Proposition \ref{th:existence of a controller psi}.

Along the solutions of (\ref{eq:main system}) with $u=\varphi_g(x_1, x_2, \tilde u)$, it yields for all $(x_1, x_2)$ in $\R^n$ and $\tilde u$ in $\R$,
\\[0.5em]$\dot V(x)\leq\varepsilon[ \alpha(M)-\alpha(V_1(x_1))]+(x_2-\varphi_1(x_1))[\tilde u $\\[0.5em]
\null\hfill$+ \Upsilon(x_1, x_2,u)]$,\\[0.5em]
where
\\[0.5em]$\Upsilon(x_1, x_2,u)=\partial_{x_1} V_1(x_1)\cdot\int_0^1\partial_{x_2} h_1(x_1,\eta_{x_1,x_2}(s),u)\,ds$\\[0.5em]
\null\hfill$+k h_2(x_1, x_2, u)- k  \partial_{x_1}\varphi_1(x_1) h_1(x_1, x_2, u)$.\\[0.5em]
With Item $2$, $3$ and $4$ of Assumption \ref{assume:bounded h}, the function $\Upsilon$ satisfies $|\Upsilon(x_1, x_2, u)| \leq \Delta(x_1, x_2)$ with% $\Delta$ given by \eqref{eq:Delta}
\\[0.5em]$\Delta(x_1,x_2)=|\partial_{x_1}V_1(x_1)|\int_0^1\Psi(x_1,\eta_{x_1,x_2}(s))\,ds$\hfill\refstepcounter{equation}$(\theequation)$\label{eq:Delta}\\[0.5em]
\null\hfill$+\Psi(x_1,x_2)k(1+|\partial_{x_1}\varphi_1(x_1)|)$\\[0.5em]
Using a particular case of the Cauchy-Schwartz inequality (i.e. $\alpha\leq \frac{1}{c}+ \frac{c}{4}\alpha^2$),
we get, for all $c>0$ %and for all $(x_1,x_2,u)$ in $\R^{n-1}\times \R\times \R$,
\\[0.5em]$(x_2-\varphi_1(x_1))\Upsilon(x_1, x_2, u) \leq \frac{1}{c}$\\[0.5em]\null\hfill$+\frac{c}{4} (x_2-\varphi_1(x_1))^2\Delta(x_1,x_2)^2.$\\[0.5em]
Consequently, it implies, that by taking %$\tilde u$ given by \eqref{eq:tildeu}
\begin{equation}\label{eq:tildeu}
\tilde u = (x_2-\varphi_1(x_1))\left[-c-\frac{c}{4} \Delta(x_1,x_2)^2\right]\ ,
\end{equation}
%\begin{equation}\label{eq:tildeu}
%\tilde u = (x_2-\varphi_1(x_1))\left[-c-\frac{c}{4} \Delta(x_1,x_2)^2\right]\ ,
%\end{equation}
it yields along the solutions of 
\begin{equation}
\label{24/10cl}
\dot x =f(x,\varphi_g(x_1, x_2, \tilde u))\ .
\end{equation}
 and for all $(x_1, x_2)$ in $\R^n$,
\\[0.5em]$
\dot V(x) \leq \varepsilon[ \alpha(M)-\alpha(V_1(x_1))]+\frac{1}{c} - c(x_2-\varphi_1(x_1))^2\ .
$\hfill\refstepcounter{equation}$(\theequation)$\label{24/10}\\[0.5em]
Note that for all $c\geq 1$, it gives,
$$
\dot V(x) \leq \varepsilon[ \alpha(M)-\alpha(V_1(x_1))]+1 - (x_2-\varphi_1(x_1))^2\ .
$$
%System (\ref{eq:main system}) in closed-loop with $u=\varphi_g(x_1, x_2, \tilde u)$,
 %where $\tilde u$ is defined in \eqref{eq:tildeu} 
% will be denoted by
%\begin{equation}
%\label{24/10cl}
%\dot x =f(x,\varphi_g(x_1, x_2, \tilde u))\ .
%\end{equation}
The function $V_1$ being proper, the set $\A_1\subset\R^n$ defined by
\\[0.5em]
$
\A_1 = \left\{x, \varepsilon\alpha(V_1(x_1))+(x_2-\varphi_1(x_1))^2
\leq \varepsilon \alpha(M)+1\right\},
$\\[0.5em]
is compact.
Moreover, selecting $c>1$,  we get, along the solutions of (\ref{24/10cl}),  $\dot V(x)<0$, for all $x$ such that $V(x)\geq\zeta$, where $\zeta$ is the positive value defined as $\zeta = \max_{x \in \A_1}\{V(x)\}$. Consequently, for all $c>1$, the set $\{x, V(x)\leq \zeta\}$ is globally asymptotically stable for (\ref{24/10cl}).

The function $\alpha$ being locally Lipschitz, we can define $K_\alpha$ its Lipschitz constant in the compact set $\{x, V(x)\leq\zeta\}$.
Hence, for all $x$ in $\{x, V(x)\leq \zeta\}$, it yields,
\\[0.5em]
$
\left|\alpha(V_1(x_1))-\alpha(V(x))\right|\leq \frac{k K_\alpha}{2} (x_2-\varphi_1(x_1))^2\ .
$\\[0.5em]
Consequently, with (\ref{24/10}) and $c>1$, we get along the solutions of (\ref{24/10cl}), for all $x$ such that $V(x)\leq\zeta$,
\\[0.5em]$
\dot V(x) \leq \varepsilon[ \alpha(M)-\alpha(V(x))]+\frac{1}{c} $\\[0.5em]
\null\hfill$- \left(c-\varepsilon\frac{k K_\alpha}{2}\right)(x_2-\varphi_1(x_1))^2\ .$\\[0.5em]
Finally, taking $c> c_g$ where
\\[0.5em]
$c_g=\max\left \{\frac{1}{\varepsilon[\alpha(M+\tilde a)-\alpha(M)]}, \varepsilon\frac{k K_\alpha}{2},1\right\}, $
\\[0.5em]
it gives, along the trajectories of (\ref{24/10cl}), for all $x$ such that $V(x)\leq\zeta$, $\dot V(x) \leq \varepsilon\left[ \alpha(M+ \tilde  a)-\alpha(V(x))\right]$.

Therefore, with $c>c_g$, for all $x$ such that $\zeta \geq V(x)>M+ \tilde a$, we get along the solutions of (\ref{24/10cl}),  $\dot V(x)<0$.
Since $c_g>1$ the same control gives also $\dot V(x)<0$ for all $x$ such that $V(x)\geq \zeta$.
Therefore the set $\{x ,\; V(x)\leq M+\tilde a\}$ in an attractor for system (\ref{eq:main system}) in closed-loop with $u=\varphi_g(x_1,x_2,\tilde u)$. Consequently, with (\ref{inclus}), the set  $\A+a\B$ contains a set that is globally and asymptotically stabilizable with the control law $\varphi_g(x_1, x_2)=\varphi_g(x_1, x_2, \tilde u)$ where $\tilde u$ is defined in (\ref{eq:tildeu}) and $c>c_g$.
This concludes the proof of Proposition \ref{th:existence of a controller psi}.
\end{proof}

\subsection{Proof of Theorem \ref{main:result}}
%\end{comment}

\begin{proof}
  Since Assumption \ref{assume:bounded h} holds, Proposition \ref{th:existence of a controller psi} applies. Let us choose the positive real number $0<a$ such that
\begin{equation}
\label{choiceof tilde c}
\max _{x\in \A+a\B} V_\ell (x)< \tilde{v}_\ell\ .
\end{equation}
Such values exist since Assumption \ref{assume:Covering} holds, and since $V_\ell$ is a proper function.

Let us consider the controller $\varphi_g$ given by Proposition \ref{th:existence of a controller psi} with this value of $a$.

Let us design a hybrid feedback law $\HF$ defining it as in Theorem \ref{main:result}, i.e., building an hysteresis of $\varphi_\ell$ and $\varphi_g$ on appropriate domains (see also \cite{GoebelSanfeliceTeel:2009} or \cite{Prieur:2001} for similar concepts applied to different control problems).

Consider an initial condition $(x(0,0),q(0,0))$ in $\R^n\times Q$, and a maximal solution $(x,q)$ of  (\ref{eq:main system}) in closed-loop with the hybrid feedback law $\HF=(Q,(C_q,D_q,\varphi_q)_{q=1,2})$. Let us assume, for the time-being, the following
\begin{lem}\label{existence hybrid time}
There exists a hybrid time $(\bar t,\bar j)$ in $\dom(x,q)$ such that $q(\bar t,\bar j)=1$ and $x(\bar t,\bar j)$ in $C_1$.
\end{lem}

Now, recalling (\ref{choiceof tilde c}) and using Assumption \ref{assume:linearization},
the sets $C_1$ is forward invariant for system (\ref{eq:main system}) in closed-loop with
$\varphi_\ell$. Thus with Lemma \ref{existence hybrid time}, we get that (\ref{eq:main system}) in closed-loop with the hybrid feedback law $\HF$ is globally asymptotically stable (since system (\ref{eq:main system}) in closed-loop with $\varphi_\ell$ is locally asymptotically stable).

Therefore to conclude the proof of Theorem \ref{main:result}, it remains to prove
Lemma \ref{existence hybrid time}. Let us prove this result by assuming the converse and exhibiting a contradiction. More precisely, let us assume that, for all $(t,j)$ in $\dom(x,q)$,
\begin{equation}
\label{contradiction}
x(t,j)\not\in C_1 \; \mbox{ or }\; q(\bar t,\bar j)=2\ .
\end{equation}
Thus, due to the expression of $D_2$, for all $(t,j)$ in $\dom(x,q)$, we have
\begin{equation}
\label{contradiction2}
x(t,j)\in D_2\setminus C_1 \; \mbox{ or } \;q(\bar t,\bar j)=2\ .
\end{equation}

If there is a time such that $x(\bar t,\bar j)\in D_2\setminus C_1$ and $q(\bar t,\bar j)=1$, then a jump occurs for the $q$-variable and, due to the expression of $G_1$, $x(\bar t,\bar j+1)\in C_1$ and $q(\bar t,\bar j+1)=2$, which is a contradiction with (\ref{contradiction}). Therefore, if $x(\bar t,\bar j)\in D_2\setminus C_1$, then $q(\bar t,\bar j)=2$. Thus we get with (\ref{contradiction2}), for all $(t,j)$ in $\dom(x,q)$,
$
x(t,j)\in D_2 \;\mbox{and }\; q(\bar t,\bar j)=2 .
$
Therefore the $x$-component is a solution of (\ref{eq:main system}) in closed-loop with $\varphi_g$ with does not enter  $C_1$. Since, with (\ref{choiceof tilde c}), $C_1$ strictly contains the set $\A$, we get the existence of a solution  of (\ref{eq:main system}) in closed-loop with $\varphi_g$ which does not converge to $\A+a\B$. This is a contradiction with the choice of the controller $\varphi_g$ satisfying the conclusion of Proposition \ref{th:existence of a controller psi}.

 This concludes the proof of Theorem \ref{main:result}.
  \end{proof}

\section{Conclusion}\label{conclusion}

A new design method has been suggested in this paper to combine a backstepping controller with a local feedback law. The class of designed controllers lies in the set of hybrid feedback laws. It
allows us to define a stabilizing control law for nonlinear control systems for which there exist some structural obstacles to the existence of classical continuous stabilizing feedback laws. More precisely, it is studied systems for which the backstepping tool for the design of stabilizers can not be applied.
\begin{comment}
 To apply this design method, three assumptions need to hold.
The first one is the existence of a controller which makes the origin of the closed-loop system locally asymptotically stable. The second assumption contains estimations on the terms which implies the fail of a controller designed using classical backstepping method. This assumption allows us to design a feedback law such that a given set is globally practically stable. The last assumption is an inclusion of the attractor into the basin of attraction of the nonlinear control system in closed loop with the local controller. Applying the main result of this paper to an example asks to check the three assumptions successively. In a future work, we will try to combine the first and the last assumption into one unique assumption which will be easier to check. %(as explained in Footnote \ref{note:1}).
\end{comment}
 % \cleardoublepage
 % \phantomsection
 % \addcontentsline{toc}{section}{Bibliography}
 % \renewcommand*{\refname}{}
  \bibliographystyle{amsplain}
  \bibliography{bibliography_2011-02-15}
\end{document}